\documentclass[12pt,reqno,a4wide]{article}
\usepackage{eurosym}
\usepackage{a4wide}
\usepackage{comment}
\usepackage{mathptmx}
\usepackage{amssymb}
\usepackage{amsfonts}
\usepackage{amsthm}
\usepackage{amsmath}
\usepackage{dsfont}
\usepackage{graphicx}
\usepackage{shadow}
\usepackage[all]{xy}
\usepackage{enumerate}
\usepackage{makeidx}
\usepackage{mathrsfs}
\usepackage{upgreek}
\usepackage{stmaryrd}
\usepackage[utf8]{inputenc}
\usepackage{array}
\usepackage{tikz-cd}
\usepackage{thumbpdf}
\usepackage[colorlinks=true,pagebackref]{hyperref}
\usepackage{textcomp}
\usepackage{chapterbib}
\usepackage{graphics}
\usepackage{longtable}
\usepackage{epsf}
\usepackage{bm}
\usepackage{adjustbox}
\usepackage{centernot}

\makeindex
\newtheorem{theorem}{Theorem}[section]
\newtheorem{lem}[theorem]{Lemma}
\newtheorem{thm}[theorem]{Theorem}
\newtheorem{prop}[theorem]{Proposition}
\newtheorem{cor}[theorem]{Corollary}
\theoremstyle{definition}
\newtheorem*{Beweis}{Proof}

\newtheorem{definition}[theorem]{Definition}
\newtheorem{defns}[theorem]{Definitions}
\newtheorem{rem}[theorem]{Remark}

\newtheorem{punto}[theorem]{}

\theoremstyle{remark}

\newtheorem{ex}[theorem]{Example}
\newtheorem{exs}[theorem]{Examples}

\input{tcilatex}
\begin{document}

\title{On $k$-Noetherian and $k$-Artinian Semirings\thanks{%
MSC2010: Primary 16Y60; Secondary 16P40, 16P20 \newline
Key Words: Semirings; Semimodules; Injective Semimodules; Noetherian
Semirings; Artinian Semirings \newline
The authors would like to acknowledge the support provided by the Deanship
of Scientific Research (DSR) at King Fahd University of Petroleum $\&$
Minerals (KFUPM) for funding this work through projects No. RG1304-1 $\&$
RG1304-2}}
\author{$%
\begin{array}{ccc}
\text{Jawad Abuhlail}\thanks{\text{Corresponding Author}} &  & \text{Rangga
Ganzar Noegraha}\thanks{\text{The paper is extracted from his Ph.D.
dissertation under the supervision of Prof. Jawad Abuhlail.}} \\
\text{abuhlail@kfupm.edu.sa} &  & \text{rangga.gn@universitaspertamina.ac.id}
\\
\text{Department of Mathematics and Statistics} &  & \text{Universitas
Pertamina} \\
\text{King Fahd University of Petroleum $\&$ Minerals} &  & \text{Jl. Teuku
Nyak Arief} \\
\text{31261 Dhahran, KSA} &  & \text{Jakarta 12220, Indonesia}%
\end{array}%
$}
\date{\today }
\maketitle

\begin{abstract}
We investigate left $k$-Noetherian and left $k$-Artinian semirings. We
characterize such semirings using $i$-injective semimodules. We prove in
particular, a partial version of the celebrated \emph{Bass-Papp Theorem} for
semiring. We illustrate our main results by examples and counter examples.
\end{abstract}

%\addcontentsline{toc}{section}{\protect\numberline{}{Introduction}}

\section*{Introduction}

\emph{Semirings} are, roughly, rings not necessarily with subtraction, and
generalize both rings and distributive bounded lattices. Semirings,{\ and
their \emph{semimodules} (defined, roughly, as modules not necessarily with
subtraction), have many important applications in several aspects of
Computer Science and Mathematics, e.g., Automata Theory \cite{HW1998},
Tropical Geometry \cite{Gla2002} and Idempotent Analysis \cite{LM2005}. Our
main reference for semirings and their semimodules is Golan's book \cite%
{Gol1999} and for rings and modules Wisbauer's book \cite{Wis1991}.}

\bigskip

Left (right) Noetherian rings, whose lattices of left (right)\ ideals
satisfy the Ascending Chain Condition,\ are well studied due to the role it
plays in simplifying the ideal structure such rings. On the other hand, left
(right)\ Artinian rings, whose lattices of left (right) ideals satisfy the
Descending Chain Condition, generalize simultaneously finite rings and rings
that are finite-dimensional vector spaces over fields. Several properties of
left (right) modules are valid only over rings with the ACC or the DCC. Some
of these properties characterize such rings, \emph{e.g.} the closure of the
class of left (right) injective modules under arbitrary direct sums
characterizes left (right) Noetherian rings \cite[3.39]{Rot2009}, and a ring
$R$ is left (right) Artinian if and only if every finitely generated left
(right) $R$-module is finitely cogenerated \cite[31.4]{Wis1991}.

\bigskip

A left (right) ideal $I$ of a semiring $S$ is called a $k$-ideal, iff $I\leq
_{S}S$ is \emph{subtractive} \cite{Hen1958} (equivalently, $I=Ker(S\overset{%
\pi _{I}}{\longrightarrow }S/I),$ where $\pi _{I}$ is the canonical
projection). In this paper, we consider the so called \emph{left }$k$\emph{%
-Noetherian semirings} (\emph{left }$k$-Artinian\emph{\ semirings}), whose
lattice of \emph{subtractive} left ideals satisfies the ACC (DCC). We
generalize several results known for left Noetherian (left Artinian)\ rings
to left $k$-Noetherian (left $k$-Artinian) semirings.

\bigskip

The paper is divided into two sections.

\bigskip

In Section 1, we collect the basic definitions, examples and preliminaries
used in this paper. In particular, we recall the definitions and basic
properties of \emph{exact sequences }introduced by the first author Abuhlail
\cite{Abu2014}.

\bigskip

In Section 2, we investigate \emph{left }$k$\emph{-Noetherian} (resp., \emph{%
left }$k$\emph{-Artinian}) \emph{semirings,} i.e. semirings satisfying the
ACC (resp., the DCC) on left $k$-ideals. In Example \ref{not-N-A}, we show
that $S:=M_{2}(\mathbb{R}^{+})$ is left $k$-Noetherian but \emph{not} left
Noetherian, and is left $k$-Artinian but \emph{not} left Artinian. In
Theorem \ref{dss-sArt}, we show that if every subtractive left ideal of a
semiring $S$ is a direct summand, then $S$ is left $k$-Artinian and left $k$%
-Noetherian. In Theorem \ref{ds-Noetherian}, we provide a partial version of
the celebrated Bass-Papp Theorem for semirings: we show that if $S$ is a
semiring with enough left $S$-$i$-injective semimodules and every direct sum
of $S$-$i$-injective left $S$-semimodules is $S$-$i$-injective, then $S$ is
left $k$-Noetherian.

\section{Preliminaries}

\label{prelim}

\qquad In this section, we provide the basic definitions and preliminaries
used in this work. Any notions from the theory of semirings and semimodules
that are not defined here can be found in {our main reference \cite{Gol1999}%
. We refer to \cite{Wis1991} for the foundations of the theory of module and
rings.}

\begin{definition}
(\cite{Gol1999}) A \textbf{semiring}%
\index{Semiring} is a datum $(S,+,0,\cdot ,1)$ consisting of a commutative
monoid $(S,+,0)$ and a monoid $(S,\cdot ,1)$ such that $0\neq 1$ and%
\begin{eqnarray*}
a\cdot 0 &=&0=0\cdot a%
\text{ for all }a\in S; \\
a(b+c) &=&ab+ac\text{ and }(a+b)c=ac+bc\text{ for all }a,b,c\in S.
\end{eqnarray*}
\end{definition}

\begin{punto}
\cite{Gol1999} Let $S$ and $T$ be semirings. The categories $_{S}\mathbf{SM}$
of \textbf{left} $S$-\textbf{semimodules} with arrows the $S$-linear maps, $%
\mathbf{SM}_{T}$ of right $S$-semimodules with arrows the $T$-linear maps,
and $_{S}\mathbf{SM}_{T}$ of $(S,T)$-bisemimodules are defined in the usual
way (as for modules and bimodules over rings). For a left $S$-semimodule $M,$
we write $L\leq _{S}M$ to indicate that $L$ is an $S$\textbf{-subsemimodule}
of $M.$
\end{punto}

\begin{defns}
\label{def-semiring}(\cite{Gol1999}) Let $(S,+,0,\cdot ,1)$ be a semiring.

\begin{itemize}
\item If the monoid $(S,\cdot ,1)$ is commutative, we say that $S$ is a
\textbf{commutative semiring.}

\item We say that the semiring $S$ is \textbf{additively idempotent}, iff $%
s+s=s$ for every $s\in S.$

\item The set of \emph{cancellative elements of }a left $S$-semimodules $M$
is defined as%
\begin{equation*}
K^{+}(M)=\{x\in M\mid x+y=x+z\Longrightarrow y=z\text{ for any }y,z\in M\}.
\end{equation*}%
We say that $M$ is a \emph{cancellative semimodule}, iff $K^{+}(M)=M.$
\end{itemize}
\end{defns}

\begin{exs}
({\cite{Gol1999})}

\begin{itemize}
\item Every ring is a cancellative semiring.

\item Any \emph{distributive bounded lattice} $\mathcal{L}=(L,\vee ,1,\wedge
,0)$ is a additively idempotent commutative semiring.

\item The set $(\mathbb{Z}^{+},+,0,\cdot ,1)$ (resp. $(\mathbb{Q}%
^{+},+,0,\cdot ,1),$ $(\mathbb{Q}^{+},+,0,\cdot ,1)$) of non-negative
integers (resp. non-negative rational numbers, non-negative real numbers) is
a cancellative commutative semiring which is not a ring.

\item $M_{n}(S),$ the set of all $n\times n$ matrices over a semiring $S,$
is a semiring.

\item $\mathbb{B}:=\{0,1\}$ with $1+1=1,$ is a an additively idempotent
commutative semiring called the \textbf{Boolean semiring}.

\item The \emph{max-plus algebra} $\mathbb{R}_{\max ,+}:=(\mathbb{R}\cup
\{-\infty \},\max ,-\infty ,+,0)$ is a semiring.

\item The \emph{log algebra} $(\mathbb{R}\cup \{-\infty ,\infty \},\oplus
,\infty ,+,0)$ is a semiring, where%
\begin{equation*}
x\oplus y=-ln(e^{-x}+e^{-y})
\end{equation*}
\end{itemize}
\end{exs}

\begin{ex}
\label{B(n,i)}(\cite[Example 1.8]{Gol1999}, \cite{AA1994}) Consider
\begin{equation*}
B(n,i):=(B(n,i),\oplus ,0,\odot ,1),
\end{equation*}%
where $B(n,i)=\{0,1,2,\cdots ,n-1\}$ and%
%($n\geq 2$) is considered
%with the binary operations $\oplus $ and $\odot $ are defined as follows:

$a\oplus b=a+b$ if $a+b<n$; otherwise, $a\oplus b=c$ is the unique natural
number $i\leq c<n$ satisfying $c\equiv a+b$ mod $(n-i);$

$a\odot b=ab$ if $ab<n$; otherwise, $a\odot b=c$ is the unique natural
number $i\leq c<n$ with $c\equiv ab$ mod $(n-i).$

Then $B(n,i)$ is a semiring. Notice that $B(n,0)=\mathbb{Z}_{n}$ (a group)
and that $B(2,1)=\mathbb{B}$ (the Boolean Algebra).
\end{ex}

\begin{ex}
(\cite[page 150, 154]{Gol1999}) Let $S$ be a semiring, $M$ be a left $S$%
-semimodule and $L\leq _{S}M.$ The \textbf{subtractive closure }of $L$ is
defined as%
\begin{equation}
\overline{L}:=\{m\in M\mid \text{ }m+\ell =\ell ^{\prime }\text{ for some }%
\ell ,\ell ^{\prime }\in L\}.  \label{L-s-closure}
\end{equation}%
One can easily check that $\overline{L}=Ker(M\overset{\pi _{L}}{%
\longrightarrow }M/L),$ where $\pi _{L}$ is the canonical projection. We say
that $L$ is \textbf{subtractive} (or a $k$\textbf{-subsemimodule}), iff $L=%
\overline{L}.$ The left $S$-semimodule $M$ is a \textbf{subtractive
semimodule}, iff every $S$-subsemimodule $L\leq _{S}M$ is subtractive.
\end{ex}

\begin{definition}
\cite[page 71]{Gol1999} Let $S$ be a semiring. A subtractive left (right)
ideal of $S$ is called a \emph{left} (\emph{right}) $k$\emph{-ideal }\cite%
{Hen1958}. We say that $S$ is a \textbf{left subtractive }(\textbf{right
subtractive})\textbf{\ semiring}, iff every left (right) ideal of $S$ is
subtractive. We say that $S$ is a subtractive semiring, iff $S$ is both left
and right subtractive.
\end{definition}

\begin{rem}
Whether a left subtractive semiring is necessarily right subtractive was an
open problem till a counterexample was given in \cite[Fact 2.1]{KNT2011}.
\end{rem}

Following \cite{BHJK2001}, we use the following definitions.

\begin{punto}
\label{variety}(cf., \cite{AHS2004})\ The category $_{S}\mathbf{SM}$ of left
semimodules over a semiring $S$ is a \emph{variety} in the sense of
Universal Algebra (closed under homomorphic images, subobjects and arbitrary
products). Whence $_{S}\mathbf{SM}$ is complete, i.e. has all limits (e.g.,
direct products, equalizers, kernels, pullbacks, inverse limits) and
cocomplete, i.e. has all colimits (e.g., direct coproducts, coequalizers,
cokernels, pushouts, direct colimits).
\end{punto}

\begin{punto}
An $S$-semimodule $N$ is a \textbf{direct summand} of an $S$-semimodule $M$ (%
\emph{i.e.} $M=N\oplus N^{\prime }$ for some $S$-subsemimodule $N^{\prime }$
of $M$) if and only if there exists $\alpha \in \mathrm{Comp}(\mathrm{End}%
(M_{S}))$ s.t. $\alpha (M)=N$ where for any semiring $T$ we set%
\begin{equation*}
\mathrm{Comp}(T):=\{t\in T\mid \text{ }\exists \text{ }\widetilde{t}\in T%
\text{ with }t+\widetilde{t}=1_{T}\text{ and }t\widetilde{t}=0_{T}=%
\widetilde{t}t\}.
\end{equation*}%
Indeed, every direct summand of $M$ is a retract of $M;$ the converse is not
true in general. Golan \cite[Proposition 16.6]{Gol1999} provided
characterizations of direct summands.
\end{punto}

\subsection*{Exact Sequences}

\bigskip

Throughout, $(S,+,0,\cdot ,1)$ is a semiring and, unless otherwise
explicitly mentioned, an $S$-module is a \emph{left }$S$-semimodule.

\bigskip

\begin{definition}
A morphism of left $S$-semimodules $f:L\rightarrow M$ is

$k$-\textbf{normal}, iff whenever $f(m)=f(m^{\prime })$ for some $%
m,m^{\prime }\in M,$ we have $m+k=m^{\prime }+k^{\prime }$ for some $%
k,k^{\prime }\in Ker(f);$

$i$-\textbf{normal}, iff $im(f)=\overline{f(L)}$ ($:=\{m\in M|\text{ }m+\ell
\in L\text{ for some }\ell \in L\}$).

\textbf{normal}, iff $f$ is $k$-normal and $i$-normal.
\end{definition}

There are several notions of exactness for sequences of semimodules. In this
paper, we use the relatively new notion of exactness introduced by Abuhlail
\cite[2.4]{Abu2014} which is stronger than that in the sense of \cite%
{Tak1982a}.

\begin{definition}
\label{Abu-exs}(\cite[2.4]{Abu2014}) A sequence
\begin{equation}
L\overset{f}{\longrightarrow }M\overset{g}{\longrightarrow }N  \label{LMN}
\end{equation}%
of left $S$-semimodules is \textbf{exact}, iff $f(L)=Ker(g)$ and $g$ is $k$%
-normal.
\end{definition}

\begin{punto}
We call a (possibly infinite) sequence of $S$-semimodules
\begin{equation}
\cdots \rightarrow M_{i-1}\overset{f_{i-1}}{\rightarrow }M_{i}\overset{f_{i}}%
{\rightarrow }M_{i+1}\overset{f_{i+1}}{\rightarrow }M_{i+2}\rightarrow \cdots
\label{chain}
\end{equation}

\emph{chain complex,} iff $f_{j+1}\circ f_{j}=0$ for every $j;$

\emph{exact,} iff each partial sequence with three terms $M_{j}\overset{f_{j}%
}{\rightarrow }M_{j+1}\overset{f_{j+1}}{\rightarrow }M_{j+2}$ is exact.

A \textbf{short exact sequence} (or a \textbf{Takahashi extension} \cite%
{Tak1982b}) of $S$-semimodules is an exact sequence of the form%
\begin{equation*}
0\longrightarrow L\overset{f}{\longrightarrow }M\overset{g}{\longrightarrow }%
N\longrightarrow 0.
\end{equation*}
\end{punto}

\section{Noetherian and Artinian Semirings}

As before, $(S,+,0,\cdot ,1)$ is a semiring and, unless otherwise explicitly
mentioned, an $S$-semimodule is a \emph{left} $S$-semimodule.

\label{sec-NA}

\begin{definition}
A left $S$-semimodule $M$ is

\textbf{Noetherian} (resp., $k$\textbf{-Noetherian}), iff $M$ satisfies the
ACC on its $S$-subsemimodules (resp., subtractive $S$-subsemimodules).

\textbf{Artinian} (resp., $k$\textbf{-Artinian}), iff $M$ satisfies the DCC
on its $S$-subsemimodules (resp., subtractive $S$-subsemimodules).

The corresponding notions for right $S$-semimodules are defined analogously.
\end{definition}

\begin{rem}
\label{ds-subtractive}Every direct summand of an $S$-semimodule is
subtractive. Let $M$ be an $S$-semimodule and $L$ a direct summand of $M.$
Then there exists $N\leq _{S}M$ such that $M=N\oplus L.$ Let $m\in M$ and $%
\ell ,\ell ^{\prime }\in L$ be such that $m+\ell =\ell ^{\prime }$. Write $m=%
\widetilde{n}+\widetilde{\ell }$ for some $\tilde{n}\in N$ and $\widetilde{%
\ell }\in L$, whence $m+\ell =(\tilde{n}+\widetilde{\ell })+\ell =\tilde{n}+(%
\widetilde{\ell }+\ell )=\ell ^{\prime }$. Since the sum $N+L$ is direct, $%
\tilde{n}=0,$ and thus $m=\widetilde{\ell }\in L$.$\blacksquare $
\end{rem}

The following result is an easy observation; however, we highlight it as it
will be used frequently in the proofs of the main results.

\begin{lem}
\label{lemint}Let $M$ be an $S$-semimodule and $N$ a \emph{subtractive }$S$%
-subsemimodules of $M.$ If $M=L\oplus K$ for some $L\leq _{S}N$ and $K\leq
_{S}M,$ then
\begin{equation*}
N=L\oplus (K\cap N).
\end{equation*}

\begin{Beweis}
Clearly, $L+(K\cap N)\subseteq N$. Let $n\in N.$ Since $M=L+K,$ there exist $%
k\in K$ and $\ell \in L$ such that $n=\ell +k.$ Since $\ell \in N$ and $N$
is subtractive, we have $k\in N,$ whence $n\in L+(K\cap N)$. So, $N=$ $%
L+(K\cap N).$ Suppose now that $\ell +k=\ell ^{\prime }+k^{\prime }$ for
some $\ell ,\ell ^{\prime }\in L$ and $k,k^{\prime }\in K\cap N$. Since the
sum $L+K$ is direct, $\ell =\ell ^{\prime }$ and $k=k^{\prime }$.$%
\blacksquare $
\end{Beweis}
\end{lem}

\begin{ex}
Let $S:=M_{2}(\mathbb{R}^{+}).$ Consider the left ideals%
\begin{equation*}
E_{1}=\left\{ \left[ {%
\begin{array}{cc}
a & 0 \\
b & 0%
\end{array}%
}\right] |\text{ }a,b\in \mathbb{R}^{+}\right\} \text{ and }E_{2}=\left\{ %
\left[ {%
\begin{array}{cc}
0 & c \\
0 & d%
\end{array}%
}\right] |\text{ }c,d\in \mathbb{R}^{+}\right\}
\end{equation*}%
and the left ideal%
\begin{equation*}
N_{\geq 1}:=\left\{ \left[ {%
\begin{array}{cc}
a & c \\
b & d%
\end{array}%
}\right] \text{ }|\text{ }a\leq c,b\leq d,a,b,c,d\in \mathbb{R}^{+}\right\} .
\end{equation*}%
Then we have $N_{\geq 1}\cap (E_{1}\oplus E_{2})=N_{\geq 1}\cap S=N_{\geq
1}, $ while $N_{\geq 1}\cap E1=\{0\}$ and $N_{\geq 1}\cap E_{2}=E_{2}.$ So,
we have%
\begin{equation*}
N_{\geq 1}\cap (E_{1}\oplus E_{2})\neq (N_{\geq 1}\cap E1)\oplus (N_{\geq
1}\cap E_{2}).
\end{equation*}%
Notice that $N_{\geq 1}\leq _{S}S$ is not subtractive, whence the condition
that $N$ is a subtractive subsemimodule of $M$ in Lemma \ref{lemint} cannot
be dropped.$\blacksquare $
\end{ex}

\begin{definition}
Let $S$ be a semiring, $M$ be a left $S$-semimodule and $N\leq _{S}M.$ A
\emph{subtractive} left $S$-subsemimodule $L\leq _{S}M$ is a \textbf{maximal
subtractive subsemimodule} of $N$ if $L\subsetneqq N$ and if $L^{\prime }$
is a subtractive subsemimodule of $M$ with $L\subseteq L^{\prime }\subseteq
N $, then $L=L^{\prime }$ or $L^{\prime }=N$.
\end{definition}

\begin{lem}
\label{max-sub}If $M$ is a $k$-Noetherian left $S$-semimodule, then every
non-zero subsemimodule of $M$ contains a maximal subtractive $S$%
-subsemimodule.

\begin{Beweis}
Let $N\leq _{S}S$ be a non-zero subsemimodule and consider
\begin{equation*}
\mathcal{I}:=\{L\lvertneqq _{S}N|\text{ }L\text{ is a subtractive
subsemimodule of }M\}.
\end{equation*}%
Notice that $L_{0}:=\{0_{M}\}\in \mathcal{I}$. If $L_{0}$ is a maximal
subtractive subsemimodule of $N,$ then we are done. Otherwise, there exists $%
L_{1}\in \mathcal{I}$ such that $L_{0}\subsetneqq L_{1}.$ If $L_{1}$ is a
maximal subtractive subsemimodule of $M,$ we are done. Otherwise, there
exists $L_{2}\in \mathcal{I}$ such that $L_{1}\subsetneqq L_{2}.$ If no such
maximal subsemimodule of $N$ exists, we obtain a non-terminating strictly
ascending chain
\begin{equation*}
L_{0}\subsetneqq L_{1}\subsetneqq L_{2}\subsetneqq \cdots \subsetneqq
L_{k}\subsetneqq L_{k+1}\subsetneqq \cdots
\end{equation*}%
of $S$-subsemimodules of $N$ which are subtractive subsemimodules of $M$,
absurd since $M$ is $k$-Noetherian.$\blacksquare $
\end{Beweis}
\end{lem}

\begin{definition}
The semiring $S$ is\textbf{\ left Noetherian} (resp., \textbf{left }$k$%
\textbf{-Noetherian}), iff $_{S}S$ is Noetherian (resp., left $k$%
-Noetherian), equivalently every ascending chain condition of left (resp.,
subtractive left) ideals of $S$ terminates;

\textbf{left Artinian} (resp., left $k$-Artinian), iff $_{S}S$ is Artinian
(resp., left $k$-Artinian), equivalently every descending chain of left
(resp., subtractive left) ideals of $S$ terminates.

The \textbf{right }($k$-)\textbf{Noetherian} and \textbf{right (}$k$\textbf{%
-)Artinian semirings} are defined analogously. A semiring which is both left
and right ($k$-)Noetherian is called ($k$-)\textbf{Noetherian}, and a
semiring which is both left and right ($k$-)Artinian is called ($k$-)\textbf{%
Artinian}.
\end{definition}

\begin{ex}
(\cite{AD1975})\ The semiring $\mathbb{Z}^{+}$ is Noetherian but not
Artinian. Setting $I_{k}:=\{0,k,k+1,k+2,\cdots \}$ yields the strictly
descending \emph{non-terminating} chain of ideals of $\mathbb{Z}^{+}:$%
\begin{equation*}
I_{1}\supsetneqq I_{2}\supsetneqq \cdots \supsetneqq I_{k}\supsetneqq
I_{k+1}\supsetneqq \cdots ,
\end{equation*}%
i.e. $\mathbb{Z}^{+}$ is not Artinian.
\end{ex}

\begin{lem}
\label{3-ideals}The only non-trivial proper subtractive left ideals of $%
S:=M_{2}(\mathbb{R}^{+})$ are%
\begin{eqnarray*}
E_{1} &=&Span\left( \left\{ \left[ {%
\begin{array}{cc}
1 & 0 \\
0 & 0%
\end{array}%
}\right] \right\} \right) =\left\{ \left[ {%
\begin{array}{cc}
a & 0 \\
b & 0%
\end{array}%
}\right] |\text{ }a,b\in \mathbb{R}^{+}\right\} \\
E_{2} &=&Span\left\{ \left[ {%
\begin{array}{cc}
0 & 0 \\
0 & 1%
\end{array}%
}\right] \right\} =\left\{ \left[ {%
\begin{array}{cc}
0 & a \\
0 & b%
\end{array}%
}\right] |\text{ }a,b\in \mathbb{R}^{+}\right\} \\
N_{r} &=&\left\{ \left[ {%
\begin{array}{cc}
ra & a \\
rb & b%
\end{array}%
}\right] |\text{ }a,b\in \mathbb{R}^{+}\right\} ,\text{ }r\in \mathbb{R}%
^{+}\backslash \{0\}.
\end{eqnarray*}
\end{lem}

\begin{Beweis}
We prove this technical lemma is three steps.

\textbf{Step I:}$\ E_{1},$ $E_{2}$ and $N_{r}$ ($r\in \mathbb{R}%
^{+}\backslash \{0\}$)\ are subtractive left ideals of $S.$

$E_{1}\leq S$ is a left ideal: for every $a,b,c,d,p,q,r,s\in \mathbb{R}^{+}$
we have
\begin{equation*}
\left[ {%
\begin{array}{cc}
p & q \\
r & s%
\end{array}%
}\right] \left[ {%
\begin{array}{cc}
a & 0 \\
b & 0%
\end{array}%
}\right] +\left[ {%
\begin{array}{cc}
c & 0 \\
d & 0%
\end{array}%
}\right] =\left[ {%
\begin{array}{cc}
pa+qb+c & 0 \\
ra+sb+d & 0%
\end{array}%
}\right] \in E_{1}.
\end{equation*}%
Moreover, $E_{1}$ is subtractive since%
\begin{equation*}
\left[ {%
\begin{array}{cc}
p & q \\
r & s%
\end{array}%
}\right] +\left[ {%
\begin{array}{cc}
a & 0 \\
b & 0%
\end{array}%
}\right] =\left[ {%
\begin{array}{cc}
c & 0 \\
d & 0%
\end{array}%
}\right]
\end{equation*}%
implies $q=0=s$ and $\left[ {%
\begin{array}{cc}
p & q \\
r & s%
\end{array}%
}\right] \in E_{1}$. Similarly, $E_{2}$ is a subtractive left ideal of $S$.

For any nonzero $r\in \mathbb{R}^{+}$, $N_{r}$ is a left ideal since (for
all $a,b,c,d,k,\ell ,m,n\in \mathbb{R}^{+}$) we have%
\begin{equation*}
\left[ {%
\begin{array}{cc}
k & \ell  \\
m & n%
\end{array}%
}\right] \left[ {%
\begin{array}{cc}
ra & a \\
rb & b%
\end{array}%
}\right] +\left[ {%
\begin{array}{cc}
rc & c \\
rd & d%
\end{array}%
}\right] =\left[ {%
\begin{array}{cc}
r(ka+\ell b+c) & ka+\ell b+c \\
r(ma+nb+d) & ma+nb+d%
\end{array}%
}\right] \in N_{r}.
\end{equation*}%
Moreover, $N_{r}\leq S$ is subtractive since
\begin{equation*}
\left[ {%
\begin{array}{cc}
k & \ell  \\
m & n%
\end{array}%
}\right] +\left[ {%
\begin{array}{cc}
ra & a \\
rb & b%
\end{array}%
}\right] =\left[ {%
\begin{array}{cc}
rc & c \\
rd & d%
\end{array}%
}\right] ,
\end{equation*}%
whence $c=a+k/r=a+\ell ,d=b+m/r=b+n.$ So, $k=r\ell ,m=rn$, and $\left[ {%
\begin{array}{cc}
k & \ell  \\
m & n%
\end{array}%
}\right] \in N_{r}$.

\textbf{Step II:}$\ E_{1},$ $E_{2}$ and $N_{r}$ ($r\in \mathbb{R}%
^{+}\backslash \{0\}$)\ are subtractive left ideals of $S.$

Let $I$ be a subtractive left ideal of $M_{2}(\mathbb{R}^{+})$ such that $%
E_{1}\subsetneqq I$. Then there exists $\left[ {%
\begin{array}{cc}
p & q \\
r & s%
\end{array}%
}\right] \in I$ such that $q\neq 0$ or $s\neq 0$, whence $\left[ {%
\begin{array}{cc}
0 & q \\
0 & s%
\end{array}%
}\right] \in I$ as $\left[ {%
\begin{array}{cc}
p & 0 \\
r & 0%
\end{array}%
}\right] \in I$ and%
\begin{equation*}
\left[ {%
\begin{array}{cc}
p & 0 \\
r & 0%
\end{array}%
}\right] +\left[ {%
\begin{array}{cc}
0 & q \\
0 & s%
\end{array}%
}\right] =\left[ {%
\begin{array}{cc}
p & q \\
r & s%
\end{array}%
}\right] \in I
\end{equation*}%
If $q\neq 0$, then%
\begin{equation*}
\left[ {%
\begin{array}{cc}
0 & 0 \\
0 & 1%
\end{array}%
}\right] =\left[ {%
\begin{array}{cc}
0 & 0 \\
1/q & 0%
\end{array}%
}\right] \left[ {%
\begin{array}{cc}
0 & q \\
0 & s%
\end{array}%
}\right] .
\end{equation*}%
If $s\neq 0$, then
\begin{equation*}
\left[ {%
\begin{array}{cc}
0 & 0 \\
0 & 1%
\end{array}%
}\right] =\left[ {%
\begin{array}{cc}
0 & 0 \\
0 & 1/s%
\end{array}%
}\right] \left[ {%
\begin{array}{cc}
0 & q \\
0 & s%
\end{array}%
}\right] .
\end{equation*}%
Either way $\left[ {%
\begin{array}{cc}
0 & 0 \\
0 & 1%
\end{array}%
}\right] \in I$, which implies $E_{2}\subseteq I$ and $I=S$. Similarly, if $%
I $ is a subtractive left ideal of $M_{2}(\mathbb{R}^{+})$ such that $%
E_{2}\subsetneqq I$, then $I=S$.

Let $r\in \mathbb{R}^{+}\backslash \{0\}$ and $I$ be a subtractive left
ideal of $M_{2}(\mathbb{R}^{+})$ such that $N_{r}\subsetneqq I$. Then there
exists $\left[ {%
\begin{array}{cc}
k & \ell \\
m & n%
\end{array}%
}\right] \in I$ such that $k\neq r\ell $ or $m\neq rn$. Without loss of
generality, assume that $k<r\ell $. Then $k+p=r\ell $ for some $p\in \mathbb{%
R}^{+}\backslash \{0\}$. Thus $\left[ {%
\begin{array}{cc}
p & 0 \\
q & 0%
\end{array}%
}\right] \in I$ or $\left[ {%
\begin{array}{cc}
p & 0 \\
0 & q%
\end{array}%
}\right] \in I$ for some $q\in \mathbb{R}^{+}$ as
\begin{equation*}
\left[ {%
\begin{array}{cc}
p & 0 \\
q & 0%
\end{array}%
}\right] +\left[ {%
\begin{array}{cc}
k & \ell \\
m & n%
\end{array}%
}\right] =\left[ {%
\begin{array}{cc}
r\ell & \ell \\
rn & n%
\end{array}%
}\right] \in I
\end{equation*}%
or
\begin{equation*}
\left[ {%
\begin{array}{cc}
p & 0 \\
0 & q%
\end{array}%
}\right] +\left[ {%
\begin{array}{cc}
k & \ell \\
m & n%
\end{array}%
}\right] =\left[ {%
\begin{array}{cc}
r\ell & \ell \\
m & m/r%
\end{array}%
}\right] \in I.
\end{equation*}%
Thus
\begin{equation*}
\left[ {%
\begin{array}{cc}
1 & 0 \\
0 & 0%
\end{array}%
}\right] =\left[ {%
\begin{array}{cc}
1/p & 0 \\
0 & 0%
\end{array}%
}\right] \left[ {%
\begin{array}{cc}
p & 0 \\
q & 0%
\end{array}%
}\right]
\end{equation*}%
or
\begin{equation*}
\left[ {%
\begin{array}{cc}
1 & 0 \\
0 & 0%
\end{array}%
}\right] =\left[ {%
\begin{array}{cc}
1/p & 0 \\
0 & 0%
\end{array}%
}\right] \left[ {%
\begin{array}{cc}
p & 0 \\
0 & q%
\end{array}%
}\right] .
\end{equation*}%
Either way we have $\left[ {%
\begin{array}{cc}
1 & 0 \\
0 & 0%
\end{array}%
}\right] \in I$, whence $E_{1}\subsetneqq I$ and $I=S$.

\textbf{Step III:}$\ E_{1},$ $E_{2}$ and $N_{r}$ ($r\in \mathbb{R}%
^{+}\backslash \{0\}$)\ are the \emph{only} subtractive left ideals of $S.$

Let $I$ be a proper non-trivial subtractive left ideal of $S.$ Then $\left[ {%
\begin{array}{cc}
k & \ell  \\
m & n%
\end{array}%
}\right] \in I\backslash \{0\}$ for some $k,\ell ,m,n\in \mathbb{R}^{+}$. If
$k\neq 0$, then
\begin{equation*}
\left[ {%
\begin{array}{cc}
1/k & 0 \\
0 & 0%
\end{array}%
}\right] \left[ {%
\begin{array}{cc}
k & \ell  \\
m & n%
\end{array}%
}\right] =\left[ {%
\begin{array}{cc}
1 & \ell /k \\
0 & 0%
\end{array}%
}\right]
\end{equation*}%
whence $\left[ {%
\begin{array}{cc}
1 & 0 \\
0 & 0%
\end{array}%
}\right] \in I,$ or $\left[ {%
\begin{array}{cc}
k/\ell  & 1 \\
0 & 0%
\end{array}%
}\right] \in I$, and it follows that $I\in \{E_{1},N_{k/\ell },S\}$ as $I$
contains $E_{1}$ or $N_{k/\ell }$. If $\ell \neq 0$, then
\begin{equation*}
\left[ {%
\begin{array}{cc}
0 & 0 \\
1/\ell  & 0%
\end{array}%
}\right] \left[ {%
\begin{array}{cc}
k & \ell  \\
m & n%
\end{array}%
}\right] =\left[ {%
\begin{array}{cc}
0 & 0 \\
k/\ell  & 1%
\end{array}%
}\right] ,
\end{equation*}%
whence $\left[ {%
\begin{array}{cc}
0 & 0 \\
0 & 1%
\end{array}%
}\right] \in I$ or $\left[ {%
\begin{array}{cc}
0 & 0 \\
k/\ell  & \ell
\end{array}%
}\right] \in I$, and so $I\in \{E_{2},N_{k/\ell },S\}$ as $I$ contains $E_{2}
$ or $N_{k/\ell }$. If $m\neq 0$, then
\begin{equation*}
\left[ {%
\begin{array}{cc}
0 & 1/m \\
0 & 0%
\end{array}%
}\right] \left[ {%
\begin{array}{cc}
k & \ell  \\
m & n%
\end{array}%
}\right] =\left[ {%
\begin{array}{cc}
1 & n/m \\
0 & 0%
\end{array}%
}\right]
\end{equation*}%
whence $\left[ {%
\begin{array}{cc}
1 & 0 \\
0 & 0%
\end{array}%
}\right] \in I$ or $\left[ {%
\begin{array}{cc}
m/n & 1 \\
0 & 0%
\end{array}%
}\right] \in I$, and it follows that $I\in \{E_{1},N_{m/n},S\}$ as $I$
contains $E_{1}$ or $N_{m/n}$. If $n\neq 0$, then
\begin{equation*}
\left[ {%
\begin{array}{cc}
0 & 0 \\
0 & 1/n%
\end{array}%
}\right] \left[ {%
\begin{array}{cc}
k & \ell  \\
m & n%
\end{array}%
}\right] =\left[ {%
\begin{array}{cc}
0 & 0 \\
m/n & 1%
\end{array}%
}\right]
\end{equation*}%
whence $\left[ {%
\begin{array}{cc}
0 & 0 \\
0 & 1%
\end{array}%
}\right] \in I$ or $\left[ {%
\begin{array}{cc}
0 & 0 \\
m/n & \ell
\end{array}%
}\right] \in I.$ So, $I\in \{E_{2},N_{m/n},S\}$ as $I$ contains $E_{2}$ or $%
N_{m/n}$.$\blacksquare $
\end{Beweis}

We provide an example of a semiring which is left $k$-Artinian and left $k$%
-Noetherian but neither left Artinian (nor left Noetherian):

\begin{ex}
\label{not-N-A}Let $S=M_{2}(\mathbb{R}^{+})$. By Lemma \ref{3-ideals}, the
only subtractive left ideals of $S$ are $0$, $S,$ $E_{1},$ $E_{2}$ and $%
N_{r} $ ($r\in \mathbb{R}^{+}\backslash \{0\}$). Notice that for $r\neq s$,
the left ideals $N_{r},N_{s}$ are not comparable. Thus, the longest
ascending (descending) chain of subtractive left ideals of $S$ is $%
0\subsetneqq N\subsetneqq S$ ($S\varsupsetneqq N\varsupsetneqq 0$) with $%
N=E_{2}$ or $N=N_{r}$ for some $r\in \mathbb{R}^{+}.$ Whence, $S$ is left $k$%
-Artinian and left $k$-Noetherian.

On the other hand, for every $r\in \mathbb{R}^{+}$ we have a left ideal of $%
S $ given by%
\begin{equation*}
N_{\geq r}=\left\{ \left[ {%
\begin{array}{cc}
a & p \\
b & q%
\end{array}%
}\right] :p\geq ra,\text{ }q\geq rb,\text{ }a,b,p,q\in \mathbb{R}%
^{+}\right\} .
\end{equation*}%
Thus, we have an infinite strictly descending chain of left ideal that does
not terminate%
\begin{equation*}
N_{1}\supsetneqq N_{\geq 2}\supsetneqq N_{\geq 3}\supsetneqq \cdots
\supsetneqq N_{\geq m}\supsetneqq N_{\geq m+1}\supsetneqq \cdots ,
\end{equation*}%
i.e. $S$ is not $k$-Artinian. On the other hand, we have an infinite
ascending chain of left ideals that does not terminate
\begin{equation*}
N_{\geq 1}\subsetneqq N_{\geq \frac{1}{2}}\subsetneqq N_{\geq \frac{1}{3}%
}\subsetneqq \cdots \subsetneqq N_{\geq \frac{1}{m}}\subsetneqq N_{\frac{1}{%
m+1}}\subsetneqq \cdots ,
\end{equation*}%
i.e. $S$ is not $k$-Noetherian.$\blacksquare $
\end{ex}

An additional example of a $k$-Noetherian semiring that is not Noetherian
was communicated to Abuhlail by T. Nam:

\begin{ex}
The semiring $\mathbb{R}^{+}[x]$ is $k$-Noetherian but not Noetherian.
\end{ex}

\begin{Beweis}
The semiring $\mathbf{B}[x],$ where $\mathbf{B}$ is the Boolean semiring, is
not Noetherian. The surjective morphism of semirings%
\begin{equation*}
f:\mathbb{R}^{+}\longrightarrow \mathbf{B},\text{ }r\mapsto \left\{
\begin{array}{ccc}
1, &  & r\neq 0 \\
&  &  \\
0, &  & r=0%
\end{array}%
\right.
\end{equation*}%
induces a surjective morphism of semirings $\mathbb{R}^{+}[x]\longrightarrow
\mathbf{B}[x],$ whence $\mathbb{R}^{+}[x]$ is not Noetherian.$\blacksquare $
\end{Beweis}

We do not know whether $k$-Artinian semirings are $k$-Noetherian. However,
we have the following interesting result.

\begin{lem}
\label{dcc-acc}A left $S$-semimodule $M$ satisfies the ACC on direct
summands if and only if $M$ satisfies the DCC on direct summands.
\end{lem}

\begin{Beweis}
$(\Longrightarrow )$\ Assume that $M$ satisfies the ACC on direct summands.
Let%
\begin{equation}
N_{1}\supseteq N_{2}\supseteq N_{3}\supseteq \cdots \supseteq N_{i}\supseteq
N_{i+1}\supseteq \cdots  \label{N-des}
\end{equation}%
be a descending chain of direct summands of $M.$ For every $i\in \mathbb{N}$%
, there exists a direct summand $L_{i}\leq _{S}M$ such that $M=N_{i}\oplus
L_{i}.$ Since $M=N_{2}\oplus L_{2}$ and $N_{1}\supseteq N_{2},$ we have by
(taking into consideration Remark \ref{ds-subtractive}):%
\begin{equation*}
N_{1}\overset{\text{Lemma \ref{lemint}}}{=}N_{2}\oplus (N_{1}\cap L_{2})%
\text{ and }M=N_{1}\oplus L_{1}=N_{2}\oplus (N_{1}\cap L_{2})\oplus L_{1}.
\end{equation*}%
Set $K_{1}:=L_{1}$ and $K_{2}:=(N_{1}\cap L_{2})\oplus L_{1}$, so that $%
N_{1}\oplus K_{1}=M=N_{2}\oplus K_{2}$ and $K_{1}\subseteq K_{2}.$

Now, $N_{2}\supseteq N_{3}$ and $M=N_{3}\oplus L_{3},$ whence $N_{2}\overset{%
\text{Lemma \ref{lemint}}}{=}N_{3}\oplus (N_{2}\cap L_{3})$ and so%
\begin{equation*}
M=N_{2}\oplus K_{2}=N_{3}\oplus (N_{2}\cap L_{3})\oplus K_{2}.
\end{equation*}%
Set $K_{3}:=(N_{2}\cap L_{3})\oplus K_{2}$, so that $M=N_{3}\oplus K_{3}$
and $K_{2}\subseteq K_{3}$. Continuing this way, we obtain an ascending
chain
\begin{equation}
K_{1}\subseteq K_{2}\subseteq K_{3}\subseteq \cdots \subseteq K_{i}\subseteq
K_{i+1}\subseteq \cdots  \label{K-asc}
\end{equation}%
of \emph{direct summands} of $_{S}M.$ By our assumption, the ascending chain
(\ref{K-asc}) terminates, whence there exists $t\in \mathbb{N}$ such that $%
K_{i}=K_{t}$ for any $i\geq t$. For any $i\geq t,$ we have $N_{t}\supseteq
N_{i},$ $M=N_{i}\oplus K_{i}$ and $N_{t}\cap K_{t}=0$ and so%
\begin{equation*}
N_{t}\overset{\text{Lemma \ref{lemint}}}{=}N_{i}\oplus (N_{t}\cap
K_{i})=N_{i}\oplus (N_{t}\cap K_{t})=N_{i},
\end{equation*}%
thus the descending chain (\ref{N-des}) terminates.

$(\Longleftarrow )$\ Assume that $M$ satisfies the DCC on direct summands.
Let%
\begin{equation}
L_{1}\subseteq L_{2}\subseteq L_{3}\subseteq \cdots \subseteq L_{i}\subseteq
L_{i+1}  \label{L-asc}
\end{equation}%
be an ascending chain of direct summands of $M$. For every $i\in \mathbb{N}$%
, there exists an $S$-subsemimodule $N_{i}\leq _{S}M$ such that $%
M=L_{i}\oplus N_{i};$ in particular $M=L_{1}\oplus N_{1}.$ Since $%
L_{1}\subseteq L_{2}$ it follows (taking into consideration Remark \ref%
{ds-subtractive}) that $L_{2}\overset{\text{Lemma \ref{lemint}}}{=}%
L_{1}\oplus (L_{2}\cap N_{1}),$ whence
\begin{equation*}
M=L_{2}\oplus N_{2}=L_{1}\oplus (L_{2}\cap N_{1})\oplus N_{2}.
\end{equation*}%
Since $L_{2}\cap N_{1}\subseteq N_{1}$ it follows that%
\begin{equation*}
N_{1}\overset{\text{Lemma \ref{lemint}}}{=}L_{2}\cap N_{1})\oplus (N_{1}\cap
(L_{1}\oplus N_{2})),
\end{equation*}%
whence%
\begin{equation*}
M=L_{1}\oplus N_{1}=L_{1}\oplus (L_{2}\cap N_{1})\oplus (N_{1}\cap
(L_{1}\oplus N_{2})).
\end{equation*}%
Setting $N_{1}^{\prime }:=N_{1}$ and $N_{2}^{\prime }:=N_{1}\cap
(L_{1}\oplus N_{2})$, we have $L_{1}\oplus N_{1}^{\prime }=M=L_{2}\oplus
N_{2}^{\prime }$ where $N_{1}^{\prime }\supseteq N_{2}^{\prime }$. Since $%
M=L_{2}\oplus N_{2}^{\prime }$ and $L_{2}\subseteq L_{3},$ it follows that $%
L_{3}\overset{\text{Lemma \ref{lemint}}}{=}L_{2}\oplus (L_{3}\cap
N_{2}^{\prime }),$ whence%
\begin{equation*}
M=L_{3}\oplus N_{3}=L_{2}\oplus (L_{3}\cap N_{2}^{\prime })\oplus N_{3}.
\end{equation*}%
Since $L_{3}\cap N_{2}^{\prime }\subseteq N_{2}^{\prime },$ we have%
\begin{equation*}
N_{2}^{\prime }\overset{\text{Lemma \ref{lemint}}}{=}(L_{3}\cap
N_{2}^{\prime })\oplus (N_{2}^{\prime }\cap (L_{2}\oplus N_{3})).
\end{equation*}%
Setting $N_{3}^{\prime }:=N_{2}^{\prime }\cap (L_{2}\oplus N_{3}),$ we have $%
N_{2}^{\prime }\supseteq N_{3}^{\prime }$ and%
\begin{equation*}
M=L_{2}\oplus N_{2}^{\prime }=L_{2}\oplus (L_{3}\cap N_{2}^{\prime })\oplus
N_{3}^{\prime }=L_{3}\oplus N_{3}^{\prime }.
\end{equation*}%
Continuing this process, we obtain a descending chain%
\begin{equation}
N_{1}^{\prime }\supseteq N_{2}^{\prime }\supseteq \cdots \supseteq
N_{i}^{\prime }\supseteq N_{i+1}^{\prime }\supseteq \cdots  \label{dN}
\end{equation}%
of \emph{direct summands }of $M$ such that $M=L_{i}\oplus N_{i}^{\prime }$
for every $i\in \mathbb{N}$. By our assumption, the descending chain (\ref%
{dN}) terminates, \emph{i.e. }there exists some $k\in \mathbb{N}$ such that $%
N_{i}^{\prime }=N_{k}^{\prime }$ for every $i\geq k$.

Now, for every $i\geq k$, we have $L_{k}\subseteq L_{i},$ $M=L_{k}\oplus
N_{k}^{\prime }$ and $L_{i}\cap N_{i}^{\prime }=0$ and so%
\begin{equation*}
L_{i}\overset{\text{Lemma \ref{lemint}}}{=}L_{k}\oplus (L_{i}\cap
N_{k}^{\prime })=L_{k}\oplus (L_{i}\cap N_{i}^{\prime })=L_{k}.
\end{equation*}%
Thus the ascending chain (\ref{L-asc}) terminates.$\blacksquare $
\end{Beweis}

A ring in which every left ideal is a direct summand is left Artinian and
left Noetherian \cite[3.4 and 4.1]{Wis1991} (in fact, left semisimple). The
following result extends this fact to semirings.

\begin{thm}
\label{dss-sArt}If every subtractive left ideal of $S$ is a direct summand,
then $S$ is left $k$-Artinian and left $k$-Noetherian.
\end{thm}

\begin{Beweis}
Assume that every subtractive left ideal of $S$ is a direct summand.

\textbf{Claim I:} $S$ is left $k$-Artinian.

Suppose that%
\begin{equation}
I_{1}\supsetneqq I_{2}\supsetneqq I_{3}\supsetneqq \cdots \supsetneqq
I_{i}\supsetneqq I_{i+1}\supsetneqq \cdots  \label{I-desc}
\end{equation}%
is a strictly descending chain of left subtractive ideals of $S$ that does
not terminate. For every $k\in \mathbb{N},$ there exists, by our assumption,
some left ideal $N_{k}\leq _{S}S$ such that $S=I_{k}\oplus N_{k}.$ The left
ideals $I_{k},N_{k}$ are non-zero as the chain does not terminate, and are
subtractive by Remark \ref{ds-subtractive}.

Since $I_{1}\supseteq I_{2}$ and $S=I_{2}\oplus N_{2}$, we have
\begin{equation*}
I_{1}\overset{\text{Lemma \ref{lemint}}}{=}I_{2}\oplus (I_{1}\cap N_{2}).
\end{equation*}%
Then $J_{1}:=I_{1}\cap N_{2}$ is a subtractive left ideal of $S,$ which is
non-zero as $I_{1}\supsetneqq I_{2},$ and $I_{1}=I_{2}\oplus J_{1}.$ Since $%
I_{2}\supseteq I_{3}$ and $S=I_{3}\oplus N_{3}$, we have
\begin{equation*}
I_{2}\overset{\text{Lemma \ref{lemint}}}{=}I_{3}\oplus (I_{2}\cap N_{3}).
\end{equation*}%
Then $I_{2}\cap N_{3}$ is a subtractive left ideal of $S,$ which is non-zero
as $I_{2}\supsetneqq I_{3},$ and $I_{1}=I_{2}\oplus J_{1}=I_{3}\oplus
J_{2}\oplus J_{1}.$ Continuing this process, we obtain at the $k$th step, a
non-zero subtractive left ideal $J_{k}\leq _{S}S$ such that%
\begin{equation*}
I_{k}=I_{k+1}\oplus J_{k}\text{ and }I_{1}=I_{k+1}\oplus J_{k}\oplus \cdots
\oplus J_{1}.
\end{equation*}

Setting $J_{i}^{\prime }:=J_{1}\oplus \cdots \oplus J_{i}$ for each $i\in
\mathbb{N}$, we have $S=J_{i}^{\prime }\oplus I_{i+1}\oplus N_{1}$ whence $%
J_{i}^{\prime }$ is subtractive (by Remark \ref{ds-subtractive}). One can
easily show that $J:=\bigcup\limits_{i\in \mathbb{N}}J_{i}^{\prime }$ is
subtractive.

By our assumption, $S=J\oplus N$ for some left ideal of $N\leq _{S}S$. Thus $%
1_{S}=j+n$ for some $j\in J$ and $n\in N$. Since $j\in J_{i}^{\prime }$ for
some $i\in \mathbb{N}$, it can be written in a unique way as $%
j=j_{1}+j_{2}+...+j_{i}$ for some uniquely determined $j_{k}\in J_{k},$ $%
k=1,2,...,i$. Since $J_{i+1}^{\prime }\subseteq J$, the sum $J_{i+1}^{\prime
}+N$ is direct, whence the sum $J_{1}+J_{2}+...+J_{i}+J_{i+1}+N$ is direct.
Setting
\begin{equation*}
K:=J_{1}\oplus ...\oplus J_{i}\oplus N,
\end{equation*}%
this means that the sum $J_{i+1}+K$ is direct. For any $s_{i+1}\in
J_{i+1}\backslash \{0\},$ we have%
\begin{equation*}
s_{i+1}=s_{i+1}1_{S}=s_{i+1}(j_{1}+j_{2}+...+j_{i}+n)=s_{i+1}j_{1}+s_{i+1}j_{2}+...+s_{i+1}j_{i}+s_{i+1}n
\end{equation*}%
where $s_{i+1}j_{k}\in J_{k}$ for $k=1,2,...,i$ and $s_{i+1}n\in N.$ It
follows that $s_{i+1}\in J_{i+1}\cap K=0$, absurd since $s_{i+1}\neq 0$. So,
the descending chain (\ref{I-desc}) terminates.

\textbf{Claim II:} $S$ is left $k$-Noetherian.

Let%
\begin{equation}
I_{1}\subseteq I_{2}\subseteq I_{3}\subseteq ...\subseteq I_{i}\subseteq
I_{i+1}\subseteq \cdots  \label{I-asc}
\end{equation}%
be an ascending chain of subtractive left ideals of $S$. Since every direct
summand of $_{S}S$ is subtractive (by Remark \ref{ds-subtractive}), it
follows from the proof of \textbf{Claim I} that $_{S}S$ satisfies DCC on
direct summands, whence $_{S}S$ satisfies ACC on direct summands by Lemma %
\ref{dcc-acc}. Since (\ref{I-asc}) is an ascending chain of subtractive left
ideals of $S,$ whence of direct summands of $_{S}S$ (by our assumption), the
chain terminates.$\blacksquare $
\end{Beweis}

\begin{ex}
Let $p$ be a prime number. Every subtractive ideal of the semiring $%
S=B(p+1,p)$ is a direct summand, and $S$ is $k$-Artinian and $k$-Noetherian.
\end{ex}

\begin{Beweis}
$S$ has no non-trivial subtractive ideals, thus every subtractive left ideal
of $S$ is a direct summand. Notice that $S$ is $k$-Artinian and $k$%
-Noetherian since it has finitely many elements.$\blacksquare $
\end{Beweis}

\begin{ex}
Let $S:=\mathbb{B}^{\mathbb{N}}$ with the canonical structure of a semiring
induced by that on $\mathbb{B}$. Then $S$ has a subtractive left ideal which
is not a direct summand and $S$ is neither $k$-Artinian nor $k$-Noetherian.
\end{ex}

\begin{Beweis}
The subtractive left ideal $\bigoplus\limits_{n\in \mathbb{N}}\mathbb{B}$ is
not a direct summand. Notice that neither the ascending chain
\begin{equation*}
\mathbb{B}\times \prod\limits_{n\geq 2}\{0\}\subsetneqq \mathbb{B}^{2}\times
\prod\limits_{n\geq 3}\{0\}\subsetneqq ...\subsetneqq \mathbb{B}^{i}\times
\prod\limits_{n\geq i+1}\{0\}\subsetneqq \mathbb{B}^{i+1}\times
\prod\limits_{n\geq i+2}\{0\}\subsetneqq ...
\end{equation*}%
nor the descending chain
\begin{equation*}
\{0\}\times \prod\limits_{n\geq 2}\mathbb{B}\supsetneqq \{0\}^{2}\times
\prod\limits_{n\geq 3}\mathbb{B}\supsetneqq ...\supsetneqq \{0\}^{i}\times
\prod\limits_{n\geq i+1}\mathbb{B}\supsetneqq \{0\}^{i+1}\times
\prod\limits_{n\geq i+2}\mathbb{B}\supsetneqq ...
\end{equation*}%
terminates, thus $S$ is neither $k$-Noetherian nor $k$-Artinian.$%
\blacksquare $
\end{Beweis}

\begin{punto}
Let $I$ be a left $S$-semimodule.

For a left $S$-semimodule $M,$ we say that $I$ is

$M$\textbf{-injective} \cite[page 197]{Gol1999}, iff for every \emph{%
injective} $S$-linear map $f:L\rightarrow M$ and any $S$-linear map $%
g:L\rightarrow I,$ there exists an $S$-linear map $h:M\rightarrow I$ such
that $h\circ f=g;$%
\begin{equation*}
\xymatrix{ 0\ar[r] & L \ar[r]^{f} \ar[d]_{g} & M \ar@{.>}[ld]^{h} \\ & I }
\end{equation*}

$M$-$i$\textbf{-injective} \cite{Alt2003}, iff for every \emph{normal
monomorphism} $f:L\rightarrow M$ and any $S$-linear map $g:L\rightarrow I,$
there exists an $S$-linear map $h:M\rightarrow I$ such that $h\circ f=g;$%
\begin{equation*}
\xymatrix{ 0\ar[r] & L \ar[rrr]^{f \, (normal)} \ar[d]_{g} & & & M
\ar@{.>}[llld]^{h} \\ & I }
\end{equation*}%
$M$-$e$\textbf{-injective }\cite{AIKN2018}, iff for every short exact
sequence $0\longrightarrow L\overset{f}{\longrightarrow }M\overset{g}{%
\longrightarrow }N\longrightarrow 0$ of left $S$-semimodules, the following
sequence $0\longrightarrow Hom_{S}(N,I)\overset{(g,I)}{\longrightarrow }%
Hom_{S}(M,I)\overset{(f,I)}{\longrightarrow }Hom_{S}(L,I)\longrightarrow 0$
of commutative monoids is exact.

We say that $I$ is \textbf{injective} (resp. $i$\textbf{-injective}, $e$%
\textbf{-injective}), iff $I$ is $M$-injective (resp. $M$-$i$-injective, $M$-%
$e$-injective) for every left $S$-semimodule $M.$
\end{punto}

\subsection*{A Bass-Papp Theorem for Semirings}

The celebrated \emph{Bass-Papp Theorem} states that a ring $R$ is left
(right) Noetherian if and only if every direct sum of left (right) injective
$R$-modules is ($R$-)injective (e.g., \cite[3.39]{Rot2009}, \cite[page 407]%
{Gri2007}). This characterization was extended to semirings: $S$ by Il'in
and

An $e$-injective version of this theorem was obtained by Abuhlail et al.:

\begin{prop}
\label{ds-ring}\emph{(\cite[Theorem 3.6]{Ili2010} ,\cite[Theorem 5.5]%
{AIKN2018}) }The following are equivalent for a semiring $S:$

\begin{enumerate}
\item $S$ is a left Noetherian \emph{ring;}

\item every direct sum of injective left $S$-semimodules is injective and
the left $S$-semimodule $S/V(S)$ can be embedded in an injective left $S$%
-semimodule.

\item every direct sum of $e$-injective left $S$-semimodules is $e$%
-injective and the left $S$-semimodule $S/V(S)$ can be embedded in an $e$%
-injective left $S$-semimodule.
\end{enumerate}
\end{prop}

Notice that the assumptions in Proposition \ref{ds-ring} force the semiring $%
S$ to be a \emph{ring}. In what follows we provide a partial version of the
Bass-Papp characterization for left $k$-Noetherian semirings.

\begin{punto}
We say that a semiring $S$ \textbf{has enough }$S$-$i$\textbf{-injective
left semimodules}, iff every left $S$-semimodule can be embedded into an $S$-%
$i$-injective left $S$-semimodule.
\end{punto}

\begin{thm}
\label{ds-Noetherian}Let $S$ be a semiring with enough left $S$-$i$%
-injective left semimodules. If every direct sum of left $S$-$i$-injective $%
S $-semimodules is $S$-$i$-injective, then $S$ is left $k$-Noetherian.
\end{thm}

\begin{Beweis}
Let
\begin{equation}
L_{1}\subseteq L_{2}\subseteq L_{3}\subseteq \cdots \subseteq L_{i}\subseteq
L_{i+1}\subseteq \cdots  \label{L-N-asc}
\end{equation}%
be a chain of subtractive left ideals of $S$ and consider the left ideal $%
L=:\bigcup\limits_{n\in \mathbb{N}}L_{n}.$ It is clear that $L$ is
subtractive. By our assumption, there exists for every $n\in \mathbb{N}$ an $%
S$-$i$-injective left $S$-semimodule $J_{n}$ and an embedding $S/L_{n}%
\overset{\iota _{n}}{\hookrightarrow }J_{n}.$ Set $J:=\bigoplus\limits_{n\in
\mathbb{N}}J_{n}$ and consider the $S$-linear map%
\begin{equation*}
\varphi _{n}:S\overset{\pi _{i}}{\longrightarrow }S/L_{n}\overset{\iota _{n}}%
{\hookrightarrow }J_{n}\text{ and }\varphi :L\rightarrow J,\text{ }x\mapsto
\sum\limits_{k=1}^{\infty }\varphi _{k}(x).
\end{equation*}%
Notice that $\varphi $ is well defined as each $x\in L$ belongs to $L_{n}$
for some $n\in \mathbb{N}$ and so $\varphi _{k}(x)=0$ for all $k\geq n,$
\emph{i.e.} $\varphi (x)=\sum\limits_{k=1}^{\infty }\varphi
_{k}(x)=\sum\limits_{k=1}^{n-1}\varphi _{k}(x)$.

By our assumption, $J$ is $S$-$i$-injective and so there exists an $S$%
-linear map $\psi :S\longrightarrow J$ such that $\psi \circ \iota =\varphi
. $
\begin{equation*}
\xymatrix{0 \ar[r] & L \ar[r]^{\iota} \ar[d]_{\varphi} & S
\ar@{-->}[ld]^{\psi} \\ & J }
\end{equation*}%
Let $\psi (1_{S})=\sum t_{k}\in J.$ Then $\psi (1_{S})\in
\bigoplus\limits_{k=1}^{m-1}J_{k}$ for some $m,$ whence $\psi (x)=\psi
(x\cdot 1_{S})=x\psi (1_{S})\in \bigoplus\limits_{k=1}^{m-1}J_{k}$ for every
$x\in L.$ In particular, $\varphi _{m}(x)=(\pi _{m}\circ \phi )(x)=0$ where $%
\pi _{m}$ is the projection on $J_{m}$. Thus $x\in L_{n}$ and $L=L_{n},$
whence $L_{k}=L_{n}$ for all $k\geq n,$ i.e. the chain terminates.
Consequently, $S$ is $k$-Noetherian.$\blacksquare $
\end{Beweis}

\begin{ex}
If $S$ is an additively idempotent semiring such that every direct sum of
left $S$-$i$-injective $S$-semimodules is $S$-$i$-injective, then $S$ is
left $k$-Noetherian. This results from Theorem \ref{ds-Noetherian} and the
fact that every semimodule over an additively idempotent semiring can be
embedded into an $e$-injective (whence an $i$-injective)\ left $S$%
-semimodule \cite[4.5]{AIKN2018}.
\end{ex}

\begin{thm}
\label{left-spit-k-Noeth}If $S$ is a semiring such that every short exact
sequence of left $S$-semimodules $0\rightarrow L\rightarrow S\rightarrow
N\rightarrow 0$ is left splitting, then $S$ is a left $k$-Noetherian.
\end{thm}

\begin{Beweis}
Let%
\begin{equation*}
N_{0}\subsetneqq N_{1}\subsetneqq N_{2}\subsetneqq \cdots \subsetneqq
N_{k}\subsetneqq N_{k+1}\subsetneqq \cdots
\end{equation*}%
be a non-terminating ascending chain of subtractive left ideals of $S$.
Notice that $N:=\bigcup\limits_{i\in \mathbb{N}}N_{i}$ is a subtractive
ideal of $S$, whence (by assumption) the following short exact sequence
\begin{equation*}
0\rightarrow N\overset{\iota }{\longrightarrow }S\overset{\pi }{%
\longrightarrow }S/N\rightarrow 0
\end{equation*}%
of left $S$-semimodules is left splitting. Let $h:S\rightarrow N$ be an $S$%
-linear map such that $h\circ \iota =id_{N}$. Then $h(1_{S})\in N$, that is $%
h(1_{S})\in N_{i}$ for some $i\in \mathbb{N}$. If $x\in N_{i+1}\backslash
N_{i},$ then%
\begin{equation*}
x=(h\circ \iota )(x)=h(x)=h(x1_{S})=xh(1_{S})\in N_{i}
\end{equation*}%
a contradiction. Hence $S$ is $k$-Noetherian.$\blacksquare $
\end{Beweis}

\begin{cor}
If $S$ is a semiring such that every subtractive left ideal is $S$-$i$%
-injective, then $S$ is a left $k$-Noetherian semiring.
\end{cor}

%\addcontentsline{toc}{section}{\protect\numberline{}{Index}} \printindex

\end{document}